\documentclass[11pt]{amsart}
\usepackage[sorted-cites]{amsrefs}
\usepackage{verbatim}
\usepackage{color}
\usepackage{amscd}
\renewcommand{\bar}{\overline}

\newcommand{\rr}{\mathbb{R}}
\newfont{\fnt}{cmr10 scaled 550}

\newtheorem{theorem}{Theorem}

\newtheorem{lemma}{Lemma}
\newtheorem{cor}{Corollary}
\newtheorem{prop}{Proposition}

\newtheorem{definition}{Definition}

\theoremstyle{remark}
\newtheorem{remark}{Remark}

\font\strange=msbm10

\renewcommand{\epsilon}{\varepsilon}

\renewcommand{\Sigma}{\varSigma}
\newcommand{\R}{{{\mathchoice  {\hbox{$\textstyle{\text{\strange R}}$}}
{\hbox{$\textstyle{\text{\strange R}}$}}
{\hbox{$\scriptstyle  N\kern-0.3em  R$}}
{\hbox{$\scriptscriptstyle  R\kern-0.2em  R$}}}}}

\newcommand{\Z}{{{\mathchoice  {\hbox{$\textstyle{\text{\strange Z}}$}}
{\hbox{$\textstyle{\text{\strange Z}}$}}
{\hbox{$\scriptstyle  Z\kern-0.3em  Z$}}
{\hbox{$\scriptscriptstyle  Z\kern-0.2em  Z$}}}}}

\newcommand{\N}{{{\mathchoice  {\hbox{$\textstyle{\text{\strange N}}$}}
{\hbox{$\textstyle{\text{\strange N}}$}}
{\hbox{$\scriptstyle  N\kern-0.3em  N$}}
{\hbox{$\scriptscriptstyle  N\kern-0.2em  N$}}}}}

\renewcommand{\phi}{\varphi}
\usepackage{amsmath,amsthm,amscd}

\newcommand{\Ric}{\textrm{Ric}}

\begin{document}

\title{Eigenvalue estimate  and compactness for  closed $f$-minimal surfaces}



\thanks{The first and third authors are  partially supported by CNPq and Faperj of Brazil. The second author is supported by CNPq of Brazil.}

\author{Xu Cheng}
\address{Instituto de Matematica e Estat\'\i stica, Universidade Federal Fluminense,
Niter\'oi, RJ 24020, Brazil. email: xcheng@impa.br}
\author{Tito Mejia}
\address{Instituto de Matematica e Estat\'\i stica, Universidade Federal Fluminense,
Niter\'oi, RJ 24020, Brazil. email: tmejia.uff@gmail.com}
\author[Detang Zhou]{Detang Zhou}
\address{Instituto de Matematica e Estat\'\i stica, Universidade Federal Fluminense,
Niter\'oi, RJ 24020, Brazil. email: zhou@impa.br}


\newcommand{\M}{\mathcal M}

\begin{abstract} Let $\Omega$ be a bounded domain with convex boundary in a complete noncompact Riemannian manifold with Bakry-\'Emery Ricci curvature bounded below by a positive constant.
We prove a lower bound of the first eigenvalue of the weighted Laplacian for closed embedded $f$-minimal hypersurfaces  contained in $\Omega$. Using this estimate, we prove a compactness theorem for the space of closed embedded $f$-minimal surfaces with the uniform upper bounds of genus and diameter in a complete $3$-manifold with Bakry-\'Emery Ricci curvature bounded below by a positive constant and admitting an exhaustion by bounded domains with convex boundary.  
\end{abstract}

\maketitle
\baselineskip=18pt

\section{Introduction}

A  hypersurface  $\Sigma$  immersed in a Riemannian manifold $(M,\bar{g})$ is said to be $f$-minimal if its   mean curvature $H$ satisfies that, for any $p\in \Sigma$,  $$H=\langle \overline{\nabla} f,\nu\rangle,$$ where $\nu$ is the unit normal at $p\in \Sigma$, $f$ is a smooth function defined on $M$, and $ \overline{\nabla} f$ denotes the gradient of $f$ on $M$.  When $f$ is a constant function,   an $f$-minimal hypersurface is just a minimal hypersurface. One  nontrivial class of $f$-minimal hypersurfaces is self-shrinker for mean curvature flow in the Euclidean space $(\rr^{n+1} , g_{can})$. Recall that a self-shrinker is a hypersurface immersed in $(\rr^{n+1} , g_{can})$ satisfying that $$H=\frac12 \langle x,\nu\rangle,$$
where $x$ is the position vector in $\mathbb{R}^{n+1}$.
Hence a self-shrinker is an $f$-minimal hypersurface $\Sigma$ with  $f=\frac{|x|^2}{4}$  (see, for instance \cite{CM3} and the references therein  the work about self-shrinkers).  

In the study of $f$-minimal hypersurfaces, it is convenient  to consider the ambient space as  a smooth metric measure space $(M,\bar{g},e^{-f}d\mu)$, where $d\mu$ is the volume form of $\bar{g}$.  For  $(M,\bar{g},e^{-f}d\mu)$, an important and natural tensor  is the Bakry-\'Emery Ricci curvature 
$\overline{\textrm{Ric}}_f:=\overline{\textrm{Ric}}+\overline{\nabla}^2f$.  
There are many interesting examples of smooth metric measure spaces $(M,\bar g,e^{-f}d\mu)$ with $\overline{\text{Ric}}_f\geq k$, where constant $k$ is positive.  A  non-trivial class of  examples is  shrinking gradient Ricci soliton. It is known that after a normalization, a shrinking gradient Ricci soliton $(M,\bar{g},f)$ satisfies the equation  $
\overline{\textrm{Ric}}+\overline{\nabla}^2f=\frac12 \bar{g}$ or equivalently  $\overline{\textrm{Ric}}_f= \frac12$. 
We refer to \cite{C}, a  survey of this topic where some compact and noncompact examples are explained.  Even though the asymptotic growth of the potential function $f$ of a noncompact shrinking gradient Ricci soliton is close to that of Gaussian shrinking soliton \cite{CZ},  both geometry and topology can be quite different.  We may consider $f$-minimal hypersurfaces in a shrinking gradient Ricci soliton. For instance, a self shrinker in $\mathbb{R}^{n+1}$ can be viewed as an $f$-minimal hypersurface in a Gauss shrinking soliton $(\mathbb{R}^{n+1}, g_{can},\frac{|x|^2}{4})$.
 
There are other examples of $f$-minimal hypersurfaces.  
  Let $M$ be the hyperbolic space $\mathbb{H}^{n+1}(-1)$.  Let  $r$ denote the distance function from a fixed point $p\in M$ and  $f(x)=nar^2(x)$, where $a>0$ is a constant. Then $\overline{\textrm{Ric}}_f\ge n(2a-1)$ and   the geodesic sphere of radius $r$ centered at $p$ in $\mathbb{H}^{n+1}(-1)$  is an $f$-minimal hypersurface if it satisfies $2ar=\coth r$.

An $f$-minimal hypersurface $\Sigma$ has two aspects to view.  One is that $\Sigma$ is $f$-minimal if and only if $\Sigma$ is a critical point of the weighted volume functional $e^{-f}d\sigma$,  where $d\sigma$ is the volume element of $\Sigma$. Another one is that $\Sigma$ is $f$-minimal if and only if $\Sigma$ is minimal in the new conformal metric $\tilde{g}=e^{-\frac{2f}{n}}\bar {g}$ (see Section \ref{notation}).  $f$-minimal hypersurfaces have been studied before, even  more general stationary hypersurfaces for parametric elliptic functionals,   see for instance the work of White \cite{W} and Colding-Minicozzi \cite{CM5}.

In this paper, we will first estimate the lower bound of the  first eigenvalue of the weighted Laplacian $\Delta_f=\Delta-\langle\nabla f,\nabla\cdot\rangle$ for closed (i.e. compact and without boundary) embedded $f$-minimal hypersurfaces   in a complete metric measure space $(M,\bar{g},e^{-f}d\mu)$.  Subsequently  using the eigenvalue estimate, we study the compactness for the space of closed embedded $f$-minimal surfaces  in a complete noncompact $3$-manifold.   To explain our result,  we  recall some backgrounds and related results.

In 1983, Choi-Wang \cite{CW} estimated the  lower bound  for  the first  eigenvalue of  closed minimal hypersurfaces in a complete Riemannian manifold  with Ricci curvature bounded below by a positive constant and proved the following
\begin{theorem}\label{CW} \cite{CW}
 If $M$ is a simply connected complete Riemannian manifold with Ricci curvature bounded below by a constant $k>0$
 and
 $\Sigma$ is a closed embedded minimal hypersurface, then the first eigenvalue of the Laplacian $\Delta$ on $\Sigma$ is at least $\frac {k}2$.
  \end{theorem}

  Later, using  a covering argument,  Choi-Schoen  \cite{CS} proved  that the assumption that $M$ is simply connected is not needed.
Recently,  Du-Ma  \cite{LiM} extended  Theorem \ref{CW}  to the first eigenvalue of the weighted Laplacian $\Delta_f$ on
 a closed embedded $f$-minimal hypersurface   in a simply connected compact manifold with positive Bakry-\'Emery Ricci curvature $\overline{\textrm{Ric}}_{f}$.  Very recently, Li-Wei \cite{LW} also used the covering argument to delete the assumption that the ambient space is simply connected in the result of  Du-Ma.

Observe that a complete manifold with  Ricci curvature bounded below by a positive constant must be compact.  But this conclusion is not true for  complete manifolds with Bakry-\'Emery Ricci curvature $\overline{\textrm{Ric}}_{f}$ bounded below by a positive constant.  One example is   Gaussian shrinking soliton  $(\rr^{n+1} , g_{can}, e^{-\frac{|x|^2}4}d\mu)$   with  $\overline{\textrm{Ric}}_{f}=\frac12$.  Hence the theorems of Du-Ma and Li-Wei cannot be applied  to  self-shrinkers.

 For self-shrinkers,
 Ding-Xin \cite{DX} recently obtained a lower bound of the first eigenvalue $\lambda_1(\mathcal{L})$ of the weighted Laplacian  $\mathcal{L}=\Delta -\frac12\langle x, \nabla \cdot\rangle$ (i.e. $\Delta_f$) on  a closed $n$-dimensional embedded self-shrinker in the Euclidean space $\mathbb{R}^{n+1}$, that is,  $\lambda_1(\mathcal{L})\ge \frac 14$.

We will discuss the  lower bound of the first eigenvalue of $\Delta_f$  of  a closed embedded $f$-minimal hypersurface in the case that 
 the ambient space is complete  noncompact. Precisely, we prove the following
\begin{theorem}\label{est1-1}
Let $(M,\bar{g},e^{-f}d\mu)$ be a complete noncompact smooth metric measure space with Bakry-\'Emery Ricci curvature $\overline{\textrm{Ric}}_{f}\geq k$, where $k$ is a positive constant.  Let  $\Sigma$ be a closed  embedded $f$-minimal hypersurface  in $M$. If there is a bounded domain $D$ in $M$ with convex boundary $\partial D$ so that $\Sigma$  is contained in  $D$,
then the first eigenvalue $\lambda_{1}(\Delta_{f}) $ of the weighted Laplacian $\Delta_f$ on $\Sigma$ satisfies
\begin{equation}
\lambda_{1}(\Delta_{f})\geq\frac{k}{2}.
\end{equation}
\end{theorem}

 The  boundary $\partial D$  is called convex if,  for any $p\in \partial D$,  the second fundamental form $A$ of $\partial D$ at $p$ is nonnegative with respect to outer unit normal of $\partial D$.

A closed self-shrinker  $\Sigma^n$ in $\rr^{n+1}$ satisfies the assumption of Theorem \ref{est1-1} since  there  always exists  a ball $D$ containing $\Sigma$. Therefore Theorem \ref{est1-1}  implies the  result of Ding-Xin for self-shrinkers mentioned before.  Besides,  we give a different and hence alternative proof of their result.

\begin{remark} If $M$ is a Cartan-Hadamard manifold, all geodesic balls are convex.
If $M$ is a complete noncompact Riemannian manifold with nonnegative sectional curvature, the  work of Cheeger-Gromoll \cite{CG} asserts that  $M$ admits  an exhaustion by convex domains. 
\end{remark}

In  \cite{CW}, 
Choi-Wang used the lower bound estimate of the first eigenvalue  in Theorem \ref{CW} to obtain an upper bound of area of a simply connected closed embedded minimal surface   $\Sigma$ in a $3$-maniflod, depending on the  genus $g$ of  $\Sigma$ and the positive lower bound $k$ of Ricci curvature of $M$. Further the lower bound of the first eigenvalue and the upper bound of area
 were used by Choi-Schoen \cite{CS} to prove a smooth  compactness theorem for   the space of closed embedded minimal surfaces of genus $g$ in a closed $3$-manifold $M^3$ with positive Ricci curvature. Very recently, Li-Wei \cite{LW} proved a compactness theorem for closed embedded $f$-minimal surfaces in a compact $3$-manifold with Bakry-\'Emery Ricci curvature $\overline{\textrm{Ric}}_{f}\geq k$, where constant $k>0$.

Recently Ding-Xin \cite{DX} applied  the lower bound estimate of the first eigenvalue of the weighted Laplacian on a self-shrinker to prove a compactness theorem  for closed self-shrinkers with uniform bounds of genus and  diameter. As it was mentioned before,  a self-shrinker in $\mathbb{R}^3$ is an $f$-minimal surface in a complete noncompact $\mathbb{R}^3$ with $\overline{\textrm{Ric}}_{f}\geq \frac 12$.  Motivated by this example,  we consider the compactness for $f$-minimal surfaces in a complete noncompact manifold. We prove that

\begin{theorem}\label{compact-intro} Let $(M^{3},\bar{g},e^{-f}d\mu)$ be a   complete noncompact smooth metric measure space with $\overline{\textrm{Ric}}_{f}\geq k$, where $k$ is a positive constant.  Assume that $M$ admits an exhaustion by bounded domains with convex boundary. Then the space, denoted by $S_{D, g}$,  of  closed embedded $f$-minimal surface in $M$ with genus  at most $g$ and diameter at most $D$ is compact in the $C^m$ topology, for any $m\geq 2$. 

Namely, any sequence in $S_{D, g}$ has a subsequence that converges in the $C^m$ topology on compact subsets to a surface in $S_{D, g}$, for any $m\geq 2$.
\end{theorem}
Theorem \ref{compact-intro} implies especially the compactness theorem of Ding-Xin for self-shrinkers. 
We also prove the following  compactness theorem, which implies Theorem \ref{compact-intro}.

\begin{theorem}\label{compact-bound} Let $(M^{3},\bar{g},e^{-f}d\mu)$ be a   complete noncompact smooth metric measure space with $\overline{\textrm{Ric}}_{f}\geq k$, where $k$ is a positive constant. Given  a   bounded domain $\Omega$, let $S$ be the space of  closed embedded $f$-minimal surface in $M$ with genus  at most $g$ and contained in the closure  $\overline{\Omega}$. If  there is  a bounded domain $U$ with convex boundary so that  $\overline{\Omega}\subset U$,  then $S$ is compact in the $C^m$ topology, for any $m\geq 2$.

Namely, any sequence in $S$ has a subsequence that converges in the $C^m$ topology on compact subsets to a surface in $S$, for any $m\geq 2$.

\end{theorem}
If $M$ admits an exhaustion by bounded domains with convex boundary, such $U$ as in Theorem \ref{compact-bound} always exists.
Also, the assumption that $f$-minimal surfaces are contained in the closure of a bounded domain $\Omega$ in Theorem \ref{compact-bound}  is equivalent to  that there is a uniform upper bound of extrinsic  diameter of $f$-minimal surfaces (see  Remark \ref{extrinsic} in Section \ref{compactness1}).

It is worth of mentioning that for self-shrinkers in $\mathbb{R}^3$,  Colding-Minicozzi \cite{CM4} proved a smooth compactness theorem for complete embedded self-shrinkers with the uniform upper bound of genus and the uniform scale-invariant area growth. In \cite{CMZ}, we generalized their result to the complete embedded $f$-minimal surfaces  in a   complete noncompact smooth metric measure space with $\overline{\textrm{Ric}}_{f}\geq k$, where constant $k>0$.

From Theorems \ref{compact-intro} and \ref{compact-bound},  we immediately have the following uniform curvature estimates respectively.
\begin{cor}\label{cor1}  Let $(M^3,\bar g,e^{-f}d\mu)  $  be a complete smooth metric measure space with $\overline{\textrm{Ric}}_{f}\geq k$, where $k$ is a positive constant. Then for any integer $g$ and a positive constant $D$, there exists a constant $C$  depending only on $M$,  $g$ and $D$   such that if $\Sigma$ is a closed embedded f-minimal surface of genus $g$ and diameter at most $D$ in $M$, the norm $|A| $  of the second fundamental form of $\Sigma$ satisfies
\[\max_{x\in \Sigma}|A|\le C.
\]

\end{cor}

\begin{cor} \label{cor2}Let $(M^{3},\bar{g},e^{-f}d\mu)$ be a   complete noncompact smooth metric measure space with $\overline{\textrm{Ric}}_{f}\geq k$, where $k$ is a positive constant.  Let $\Omega$ be  a   bounded domain whose closure is contained in   a bounded domain $U$ with convex boundary. Then for any integer $g$,  there exists a constant $C$  depending only on $U$,  $g$   such that if $\Sigma$ is a closed embedded f-minimal surface of genus $g$ contained in $\bar\Omega$,  the norm $|A| $  of the second fundamental form of $\Sigma$ satisfies
\[\max_{x\in \Sigma}|A|\le C.
\]

\end{cor}

On the other hand,  observe that a similar  argument of Theorem \ref{est1-1} also works for the case that the ambient space is a compact manifold with convex boundary. Hence we have  the following estimate:
\begin{theorem}\label{estloc1}
Let $(M,\bar{g})$ be a  simply connected compact manifold with convex boundary $\partial M$,  and $f$ a smooth function on $M$.  Assume that  $\overline{\textrm{Ric}}_{f}\geq k$, where $k$ is a positive constant.  If  $\Sigma$ is a closed  $f$-minimal hypersurface embedded in $M$ and does not intersect the boundary $\partial M$, then the first eigenvalue of the weighted Laplacian on $\Sigma$ satisfies
\begin{equation}
\lambda_{1}(\Delta_{f})\geq\frac{k}{2}.
\end{equation}
\end{theorem}
A special case of Theorem \ref{estloc1} is that
\begin{cor}\label{estloc2}
Let $(M,\bar{g})$ be a  simply connected compact manifold with convex boundary $\partial M$. Assume that  $M$ has   the Ricci curvature $\overline{\textrm{Ric}}\geq k$, where $k$ is a positive constant.  If  $\Sigma$ is a closed  embedded minimal hypersurface in $M$ and does not intersect the boundary $\partial M$, then the first eigenvalue of the Laplacian on $\Sigma$ satisfies
\begin{equation*}
\lambda_{1}(\Delta)\geq\frac{k}{2}.
\end{equation*}
\end{cor}

Corollary \ref{estloc2}  can be viewed as an extension of Theorem  \ref{CW} by Choi-Wang  if we consider  the empty boundary  as a convex one.

The rest of this paper is organized as follows: In Section \ref{notation} some definitions and notation are given; In Section \ref{reilly} we give some  facts which will be used later; In Section \ref{eigenvalue}  we prove Theorems \ref{est1-1} and  \ref{estloc1}; In Section \ref{compactness1} we  prove Theorems \ref{compact-intro} and \ref{compact-bound};  In Appendix we give the proof of the known Reilly formula for the weighted metric measure space, for the sake of completeness of proof.

\section{Definitions and notation}\label{notation}

In general, a smooth metric measure space, denoted by $(N,g,e^{-w}d\text{vol})$, is a Riemannian manifold $(N,g)$ together with a weighted volume form $e^{-w}d\text{vol}$ on $N$, where $w$ is a smooth function on $N$ and $d\text{vol}$ the volume element induced by the Riemannian metric $g$. The associated weighted Laplacian ${\Delta}_{w}$ is defined  by
$$
{\Delta}_{w}u:={\Delta} u-\langle{\nabla} w,{\nabla} u\rangle
$$
where ${\Delta}$ and ${\nabla}$ are  the Laplacian and  gradient on $(N,g)$ respectively.  

The second order operator $\Delta_w$  is a self-adjoint operator on the space of square integrable functions on $N$ with respect to the measure $e^{-w}d\text{vol}$.  For a closed manifold $N$, the first eigenvalue of $\Delta_w$, denoted by $\lambda_1(\Delta_w)$, is the lowest nonzero real number $\lambda_1$ satisfying
$$\Delta_w u=-\lambda_1u, \quad \text{on}\quad N.$$ 
It is well known that the the definition of $\lambda_1(\Delta_w)$  is equivalent to that 
$$\lambda_{1}(\Delta_w)=\inf_{\int_{N}ue^{-w}d\sigma=0, u\not\equiv 0}\frac{\displaystyle\int_{N}|\nabla u|^{2}e^{-w}d\sigma}{\displaystyle\int_{N}u^{2}e^{-w}d\sigma}.$$

The $\infty$-Bakry-\'Emery Ricci curvature tensor $\text{Ric}_{w}$ (for simplicity, Bakry-\'Emery Ricci curvature)  on  $(N,g,e^{-w}d\text{vol})$ is defined by
$$
\text{Ric}_{w }:=\text{Ric}+{\nabla}^{2}w
$$
where $\text{Ric}$ denotes the Ricci curvature  of $(N,g)$ and ${\nabla}^{2}w$ is the Hessian of $w$ on $N$.
If $w$ is constant,  ${\Delta}_{w}$ and $\textrm{Ric}_{w}$ are the Laplacian $\Delta$ and Ricci curvature  $\textrm{Ric}$ on $N$ respectively.

Now,  let $(M^{n+1}, \bar{g})$ be an $(n+1)$-dimensional Riemannian manifold.  Assume that $f$ is a smooth function on $M$ so that $(M^{n+1}, \bar{g}, e^{-f}d\mu)$ is a smooth metric measure space, where $d\mu$ is the volume element induced by  $\bar{g}$.

Let $i: \Sigma^{n}\to M^{n+1}$ be an $n$-dimensional smooth immersion.
Then $i: (\Sigma^{n}; i^*\bar{g}) \to (M^{n+1}, \bar{g})$ is an isometric immersion with the induced   metric $i^*\bar{g}$. For simplicity, we still  denote $i^*\bar{g}$  by $\bar{g}$ whenever there is no confusion. Let  $d\sigma$ denote the volume element of $(\Sigma,\bar{g})$. Then the function $f$ induces a weighted measure  $e^{-f}d\sigma$ on $\Sigma$. Thus we have   an induced smooth metric measure space $(\Sigma^{n}, \bar{g}, e^{-f}d\sigma)$. 

 In this paper, unless otherwise specified, we denote by a bar all quantities on $(M, \bar{g})$, for instance by $\overline{\nabla}$ and $\overline{\textrm{Ric}}$, 
 the Levi-Civita connection  and  the  Ricci curvature tensor of $(M, \bar{g})$ respectively. Also we denote for example by ${\nabla}$,  $\textrm{Ric}$, $\Delta$ and $\Delta_f$, 
 the Levi-Civita connection, the  Ricci curvature tensor,  the Laplacian, and the weighted Laplacian on $(\Sigma, \bar{g})$ respectively.

 Let $p\in \Sigma$ and  $\nu$ a unit normal at $p$. The second fundamental form $A$, the mean curvature $H$,  and the mean curvature vector  ${\bf H}$ of hypersurface $(\Sigma,\bar{g})$ are  defined respectively by:  $$A: T_p\Sigma\to T_p\Sigma, \quad A(X):=\overline{\nabla}_X\nu, \quad X\in T_p\Sigma,$$  $$H:=\text{tr}A=-\displaystyle\sum_{i=1}^n\langle \overline{\nabla}_{e_i}e_i,\nu\rangle,$$ $${\bf H}:=-H\nu.$$  Define  the weighted mean curvature vector ${\bf H}_f$ and the weighted mean curvature $H_f$ of $(\Sigma,\bar{g})$ by $${\bf H}_f:={\bf H}-(\overline{\nabla }f)^{\perp},$$
 $${\bf H}_f=-H_f\nu,$$ where $\perp$ denotes the projection to the normal bundle of $\Sigma$. It holds that  $$H_f=H-\langle \overline{\nabla} f,\nu\rangle.$$

\begin{definition} A hypersurface $\Sigma$ immersed in $(M^{n+1}, \bar{g}, e^{-f}d\mu)$ with the induced metric $\bar{g}$  is called an $f$-minimal hypersurface if  its weighted mean curvature $H_f$ vanishes identically, or equivalently if it satisfies
\begin{equation} H=\langle \overline{\nabla} f,\nu\rangle.
\end{equation}
\end{definition}

\begin{definition}The weighted volume of $(\Sigma, \bar{g})$ is defined by
\begin{equation}V_f(\Sigma):=\int_\Sigma e^{-f}d\sigma.
\end{equation}

\end{definition}

It is well known  that  $\Sigma$ is $f$-minimal if and only if $\Sigma$ is a critical point of the weighted volume functional. Namely,  it  holds that
\begin{prop}\label{first-var}
If $T$ is a compactly supported normal variational field on $\Sigma$ (i.e. $T=T^{\perp}$), then the first variation formula of the weighted volume of $(\Sigma,\bar{g})$ is given by
\begin{equation}\label{prop-fir-var}
\frac{d}{dt}V_{f}(\Sigma_{t})\biggr|_{t=0}=-\int_{\Sigma}\langle T, {\bf H}_{f}\rangle_{\bar{g}} e^{-f}d\sigma.
\end{equation}
\end{prop}

On the other hand,  an $f$-minimal submanifold can be viewed as a minimal submanifold under a conformal metric. Precisely,  define the new metric  $\tilde{g}=e^{-\frac{2}{n}f}\bar{g}$ on $M$, which is  conformal to $\bar{g}$.  Then the immersion $i: \Sigma\to M$ induces a metric $i^*\tilde{g}$ on $\Sigma$ from $(M, \tilde{g})$.   
In the following, $i^*\tilde{g}$ is still denoted by $\tilde{g}$ for simplicity of notation. The volume of $(\Sigma, \tilde{g})$   is
\begin{equation}\label{two-vol}
\tilde V(\Sigma):=\int_\Sigma d\tilde{\sigma}=\int_\Sigma e^{-f}d\sigma=\text{V}_f(\Sigma).
\end{equation}

\noindent Hence Proposition \ref{first-var} and (\ref{two-vol}) imply that 
\begin{equation}\label{fir-variation}
\int_{\Sigma}\langle T, \tilde{{\bf H}}\rangle_{\tilde{g}} d\tilde{\sigma}=\int_{\Sigma}\langle T, {\bf H}_{f}\rangle_{\bar{g}} e^{-f}d\sigma,
\end{equation}
where $d\tilde{\sigma}=e^{-f}d\sigma$ and $ \tilde{{\bf H}}$ denote the volume element and the mean curvature vector of $\Sigma$ with respect to the conformal metric $\tilde{g}$ respectively.

(\ref{fir-variation}) implies that $ \tilde{{\bf H}}=e^{\frac{2f}{n}}{\bf H}_f$. Therefore $(\Sigma, \bar{g})$ is $f$-minimal in $(M,\bar{g})$ if and only if $(\Sigma, \tilde{g})$ is minimal in $(M,\tilde{g})$.

In this paper,  for closed hypersurfaces, we choose $\nu$ to be  the outer unit normal.

\section{Some  facts on the weighted Laplacian and $f$-minimal hypersurfaces}\label{reilly}
In this section we give some known results which will be used later in this paper.  Recall  that Reilly \cite{R} proved an integral version of the Bochner formula for compact domains of a Riemannian manifold, which is called Reilly formula.  In \cite{LiM},  Du-Ma obtained a Reilly formula for metric measure spaces, which is the following proposition.  We include its  proof in Appendix  for the sake of completeness.

\begin{prop}\label{prop1}\cite{LiM}
Let $\Omega$ be a compact Riemannian manifold with boundary   $\partial\Omega $ and   $(\Omega, \bar{g},e^{-f}d\mu)$ be  a smooth metric measure space.  Then we have
\begin{eqnarray*}
\int_{\Omega}(\bar{\Delta}_{f}u)^{2}e^{-f}&=&\int_{\Omega}|\overline{\nabla}^{2}u|^{2}e^{-f}+\int_{\Omega}\overline{\textrm{Ric}}_{f}
(\overline{\nabla}u,\overline{\nabla}u)e^{-f}+2\int_{\partial\Omega}u_{\nu}\Delta_{f}(u)e^{-f}\\
& &+\int_{\partial\Omega}A(\nabla u,\nabla u)e^{-f}+\int_{\partial\Omega}u_{\nu}^{2}H_{f}e^{-f},
\end{eqnarray*}
where $\nu$ is the outward pointing unit normal to $\partial\Omega$ and $A$ is the second fundamental form of $\partial\Omega$ with respect to the normal $\nu$,  the quantities with bar denote the ones on $(\Omega, \bar{g})$ (for instance, $\overline{\textrm{Ric}}_{f}$ denotes the Bakry-\'Emery Ricci curvature on $(\Omega, \bar{g})$),  and  $\Delta_f$ and $H_f$ denote the weighted Laplacian on $\partial\Omega$ and the weighted mean curvature of $\partial\Omega$ respectively.
\end{prop}

A Riemannian manifold with Bakry-\'Emery Ricci curvature bounded below by a positive constant has some properties similar to a Riemannian manifold with Ricci curvature bounded below by a positive constant. We refer to the work of Wei-Wylie \cite{WW}, and Munteanu-Wang \cites{MW, MW1} and the references therein. 

We will use  the following proposition by Morgan \cite{M} (see also  \cite{WW} Corollary 5.1).

\begin{prop}\label{group} If a complete smooth metric measure space  $(N, g, e^{-\omega}d\text{vol})$ has $\text{Ric}_{w}\geq k$, where $k$ is a positive constant, then $N$ has finite weighted volume and finite fundamental group.
\end{prop}

For $f$-minimal hypersurfaces, the following intersection theorem holds, which was proved by Wei-Wylie (\cite{WW} Theorem 7.4).

\begin{prop}  \cite{WW} \label{inter}Any two closed $ f$-minimal hypersurfaces in a smooth metric measure space $(M,\bar{g}, e^{-f}d\mu)$ with $\overline{\text{Ric}}_f > 0$ must intersect. Thus a closed $ f$-minimal hypersurface in $M$ must be connected.
\end{prop}

The first and third authors  \cite{XZ} of the present paper proved that the finite weighted volume of a self-shrinker  immersed in $\mathbb{R}^{m}$ implies it is properly immersed.  In \cite{CMZ}, we generalizes this result to $f$-minimal submanifolds.

\begin{prop}\cite{CMZ}\label{prop-imer1} Let $\Sigma^{n}$ be an $n$-dimensional complete $f$-minimal submanifold immersed in an $m$-dimensional  Riemannian manifold  $M^m$, $n<m$. If $\Sigma$ has finite weighted volume, then  $\Sigma$ is properly immersed in $M$.
\end{prop}
\begin{remark} We studied  $f$-minimal submanifolds in \cite{CMZ}. An $f$-minimal hypersurface is an $f$-minimal submanifold with  codimension $1$.
\end{remark}

\section{Lower bound for $\lambda_{1}(\Delta_{f})$}\label{eigenvalue}
In this section, we apply  Reilly  formula  for metric measure space to prove Theorems \ref{est1-1} and  \ref{estloc1}. 

\noindent\textbf{Proof of Theorem \ref{est1-1}}. Since $\overline{\Ric}_{f}\geq k$,  where constant $k>0$,  Proposition \ref{group} implies that  $M$ has finite fundamental group.  We first assume that $M$ is simply connected.
 Since $\Sigma$ is connected (Proposition \ref{inter}) and embedded in $M$,  $\Sigma$ is orientable and divides $M$ into two components (see its proof in \cite{CS}). Thus $\Sigma$ divides $D$ into bounded two components $\Omega_1$ and $\Omega_2$. That is $D\backslash \Sigma=\Omega_1\cup\Omega_2$ with $\partial\Omega_1=\Sigma$ and  $\partial\Omega_{2}=\partial D\cup\Sigma$.

For simplicity, we denote by  $\lambda_{1}$ the first eigenvalue $\lambda_{1}(\Delta_{f})$ of the weighted Laplacian $\Delta_{f}$ on $\Sigma$. Let $h$ be a corresponding eigenfunction so that on $\Sigma$,
\begin{equation}\label{eigenfunction}
\Delta_{f}h+\lambda_{1}h=0\quad \textrm{ with}\quad\int_{\Sigma}h^{2}e^{-f}=1.
\end{equation}
 Consider the solution of the Dirichlet problem on $\Omega_{1}$ so that
\begin{equation}\label{solution-diri}\left\{
        \begin{array}{ll}
         \bar{\Delta}_{f}u=0 & \hbox{in $\Omega_{1} $} \\
         u=h& \hbox{on $ \partial\Omega_{1}=\Sigma$. }
        \end{array}
       \right.
       \end{equation}
      Substituting $\Omega_1$ for $\Omega$ and putting the solution $u$ of (\ref{solution-diri}) in Proposition \ref{prop1}. Then the assumption of $\overline{\Ric}_{f}$ implies that
\begin{equation*}
0\geq k\int_{\Omega_{1}}|\overline{\nabla}u|^{2}e^{-f}-2\lambda_{1}\int_{\Sigma}u_{\nu}he^{-f}+\int_{\Sigma}A(\nabla h,\nabla h)e^{-f},
\end{equation*}
where $\nu$ is the outer unit normal of $\Sigma$ with respect to $\Omega_1$.
By the Stokes' theorem and (\ref{solution-diri}),
$$
\int_{\Sigma}u_{\nu}he^{-f}=\int_{\Omega_{1}}(|\overline{\nabla}u|^{2}+u\bar{\Delta}_{f}u)e^{-f}=\int_{\Omega_{1}}|\overline{\nabla}u|^{2}e^{-f}.
$$
Thus
\begin{equation*}
0\geq(k-2\lambda_{1})\displaystyle\int_{\Omega_{1}}|\overline{\nabla}u|^{2}e^{-f}+\int_{\Sigma}A(\nabla h,\nabla h)e^{-f}.
\end{equation*}
If  $\int_{\Sigma}A(\nabla h,\nabla h)e^{-f}\geq0$,  by $u\not\equiv C$, we have
\begin{equation*}
\lambda_{1}\geq\frac{k}{2}.
\end{equation*}
If $\int_{\Sigma}A(\nabla h,\nabla h)e^{-f}<0$,  we consider the compact domain $\Omega_{2}$ with the boundary $\partial\Omega_2=\Sigma\cup\partial D$. Let $u$ be the solution of   the  mixed problem
\begin{equation}\label{solution-mix}\left\{
        \begin{array}{lll}
         \bar{\Delta}_{f}u=0 & \hbox{in $\Omega_{2} $} \\
         u=h& \hbox{on $ \Sigma$}\\
         u_{\tilde{\nu}}=0 &\hbox{on  $\partial D$,}
        \end{array}
       \right.
       \end{equation}
       where $\tilde{\nu}$ denotes  the outer unit normal  of $\partial D$ with respect to $\Omega_2$.

Substituting $\Omega_2$ for $\Omega$ and putting the solution  $u$ of (\ref{solution-mix}) in Proposition \ref{prop1}, we have
\begin{align*}
0&\geq\int_{\Omega_{2}}|\overline{\nabla}^{2}u|^{2}e^{-f}+k\int_{\Omega_{2}}|\overline{\nabla}u|^{2}e^{-f}-2\lambda_{1}\int_{\Sigma}hu_{\tilde{\nu}}e^{-f}\\
& +\int_{\Sigma}\tilde{A}(\nabla h,\nabla h)e^{-f}+\int_{\partial D}\tilde{A}(\nabla u,\nabla u)e^{-f},\nonumber
\end{align*}
 where $\tilde{\nu}$ denotes  the outer unit normal  of $\Sigma$ with respect to $\Omega_2$, and $\tilde{A}$ denotes the second fundamental form of $\Sigma$ with respect to normal $\tilde{\nu}$.
 
On the other hand,   the Stokes'  theorem and (\ref{solution-mix}) imply
\begin{equation*}
\int_{\Omega_{2}}|\overline{\nabla}u|^{2}e^{-f}=\int_{\partial\Omega_{2}}uu_{\tilde{\nu}}e^{-f}=\int_{\Sigma}hu_{\tilde{\nu}}e^{-f}.
\end{equation*}
Thus, we have
\begin{equation}\label{ine-1}
0\geq(k-2\lambda_{1})\int_{\Omega_{2}}|\overline{\nabla}u|^{2}e^{-f}+\int_{\Sigma}\tilde{A}(\nabla h,\nabla h)e^{-f}+\int_{\partial D}\tilde{A}(\nabla u,\nabla u)e^{-f}.
\end{equation}
Since $\partial D$ is assumed  convex,  the last term on the right side of (\ref{ine-1}) is nonnegative.
Observe that the orientations of $\Sigma$ are opposite for $\Omega_1$ and $\Omega_2$.  Namely, $\tilde{\nu}=-\nu$. Then $\tilde{A}(\nabla u,\nabla u)=-A(\nabla u,\nabla u)$ on $\Sigma$. This implies that the second term on the right side of (\ref{ine-1}) is nonnegative. Thus
\begin{equation*}
0\geq(k-2\lambda_{1})\int_{\Omega_{2}}|\overline{\nabla}u|^{2}e^{-f}.
\end{equation*}
Since $u$ is not constant function, we conclude that $k-2\lambda_{1}\leq0$.
Again, we have
\begin{equation*}
\lambda_{1}\geq\frac{k}{2}.
\end{equation*}
Therefore we obtain that $\lambda_{1}(\Delta_f)\geq\frac{k}{2}$ if $M$ is simply connected.

Second, if $M$ is not simply connected, we consider the universal covering $\hat{M}$, which is a finite $|\pi_1|$-fold covering. $\hat{M}$ is simply connected and the covering map $\pi: \hat{M}\to M$ is a locally isometry. 

Take $\hat{f}=f\circ\pi$. Obviously $\hat{M}$ has $\hat{Ric}_{\hat{f}}\geq k$,  and  the lift $\hat{\Sigma}$ of $\Sigma$ is also $\hat{f}$-minimal, embedded and closed. By Proposition \ref{inter}, $\hat{\Sigma}$ must be connected. Since $\hat{M}$ is simply connected, the closed embedded connected $\hat{\Sigma}$ must be orientable and thus divides $\hat{M}$ into two components.
Moreover the connectedness of $\hat{\Sigma}$ implies that the lift $\hat{D}$ of $D$ is also connected. Also  $\partial \hat{D}=\widehat{\partial D}$ is  smooth and convex.  Hence the assertion obtained for the simply connected ambient space  can be applied here. Thus the first eigenvalue of the weighted Laplacian $\hat\Delta_{\hat{f}}$ on $\hat{\Sigma}$ satisfies $ \lambda_{1}(\hat\Delta_{\hat{f}})\geq\frac{k}{2}$. 

Observing the lift of the first eigenfunction of $\Sigma$ is also an eigenfunction of $\hat{M}$, we have
$$\lambda_{1}(\Delta_{f})\geq\lambda_{1}(\hat\Delta_{\hat{f}})\geq\frac{k}{2}.$$
\qed

\begin{remark} \label{C^1} In Theorem \ref{est1-1},  the boundary $\partial D$ is not necessarily smooth. $\partial D$ can be assumed to be $C^1$, which is sufficient to the existence of the solution of the mixed problem  (\ref{solution-mix}).
\end{remark}

 Theorem \ref{estloc1} holds by the same argument as that of Theorem \ref{est1-1}.

\section{Upper bound on Area and total curvature of $f$-minimal surfaces}\label{sec.est-area}
In this section, we study  surfaces in a $3$-manifold. First we estimate the corresponding  upper bounds on the area and weighted area of an embedded closed $f$-minimal surface by applying the  first eigenvalue estimate  in Section \ref{eigenvalue}.  Next we discuss the    upper bound on the total curvature.  We begin with  a result of Yang and Yau \cite{YY}:
\begin{prop}\cite{YY}\label{yang}
Let $\Sigma^{2}$ be a closed orientable Riemannian surface  with genus $g$ . Then the first eigenvalue $\lambda_1(\Delta)$ of the Laplacian $\Delta$ on $\Sigma$ satisfies that
\begin{equation*}
\lambda_{1}(\Delta)\textrm{Area}(\Sigma)\leq 8\pi(1+g).
\end{equation*}
\end{prop}
Using Theorem \ref{est1-1} and Proposition \ref{yang}, we  obtain the following area estimates for  closed embedded $f$-minimal surfaces if the ambient space is simply connected.
\begin{prop}\label{est2}
Let  $(M^{3},\bar{g},e^{-f}d\mu)$ be  a simply connected complete smooth metric measure space with $\overline{\textrm{Ric}}_{f}\geq k$, where $k$ is a positive constant.  Let $\Sigma^{2}\subset M$  be a closed embedded $f$-minimal surface with genus $g$. If $\Sigma$ is  contained in a bounded domain $D$ with  convex boundary $\partial D$,  then its area and weighted area  satisfy the following inequalities respectively.
\begin{equation}\label{area}
\textrm{Area}(\Sigma)\leq \frac{16\pi(1+g)}{k}e^{\textrm{osc}_{\Sigma}f},
\end{equation}
\begin{equation}\label{f-area}
\textrm{Area}_f(\Sigma)\leq \frac{16\pi(1+g)}{k}e^{-\textrm{inf}_{\Sigma}f},
\end{equation}
where $\textrm{osc}_{\Sigma}f=\sup_{\Sigma}f-\inf_{\Sigma}f$.
\end{prop}
\noindent\textbf{Proof}. Consider  the conformal metric  $\tilde{g}=e^{-f}\bar{g}$ on $M$.  Let $\lambda_{1}(\tilde{\Delta})$ be the first eigenvalue of the Laplacian  $\tilde{\Delta}$ on $(\Sigma,\tilde{g})$,  which satisfies
$$
\lambda_{1}(\tilde{\Delta})=\inf_{\int_{\Sigma}ud\tilde{\sigma}=0,u\not\equiv 0}\frac{\displaystyle\int_{\Sigma}|\tilde{\nabla}u|^{2}_{\tilde{g}}d\tilde{\sigma}}
{\displaystyle\int_{\Sigma}u^{2}d\tilde{\sigma}},
$$
where $\tilde{\Delta}$, $\tilde{\nabla}$ and $d\tilde{\sigma}$ are the Laplacian, gradient and area element of $\Sigma$ with respect to  the metric $\tilde{g}$ respectively.

On the other hand, the first eigenvalue of the weighted Laplacian $\lambda_{1}(\Delta_{f})$ on  $(\Sigma,\bar{g})$ satisfies
$$
\lambda_{1}(\Delta_{f})=\inf_{\int_{\Sigma}ue^{-f}d\sigma=0, u\not\equiv 0}\frac{\displaystyle\int_{\Sigma}|\nabla u|_{\bar{g}}^{2}e^{-f}d\sigma}{\displaystyle\int_{\Sigma}u^{2}e^{-f}d\sigma}.
$$
Since $\tilde{\nabla}u=e^{f}\nabla u$,  $d\tilde{\sigma}=e^{-f}d\sigma$ and ${\tilde{g}}=e^{-f}{\bar{g}}$,
\begin{align*}
\lambda_{1}(\tilde{\Delta})&=\inf_{\int_{\Sigma}ue^{-f}d\sigma=0,u\not\equiv 0}\frac{\displaystyle\int_{\Sigma}|\nabla u|_{\bar{g}}^{2}d\sigma}{\displaystyle\int_{\Sigma}u^{2}e^{-f}d\sigma}\\
&\geq\inf_{\int_{\Sigma}ue^{-f}d\sigma=0, u\not\equiv 0}\frac{\displaystyle\int_{\Sigma}|\nabla u|_{\bar{g}}^{2}e^{-f+\inf_{\Sigma}(f)}d\sigma}{\displaystyle\int_{\Sigma}u^{2}e^{-f}d\sigma}\\
&= e^{\inf_{\Sigma}f}\lambda_{1}(\Delta_{f}).
\end{align*}
Combining Theorem \ref{est1-1} with  Proposition \ref{yang} gives the following estimate
\begin{equation}
\textrm{Area}(\Sigma,\tilde{g})\leq\frac{16\pi(1+g)}{k}e^{-\inf_{\Sigma}(f)}.
\end{equation}
Since $\textrm{Area}_f(\Sigma)=\int_{\Sigma}e^{-f}d\sigma=\textrm{Area}(\Sigma,\tilde{g})$,
\begin{equation*}
\textrm{Area}_f(\Sigma)\leq\frac{16\pi(1+g)}{k}e^{-\inf_{\Sigma}(f)},
\end{equation*}
which is  (\ref{f-area}).
Thus
$$
\textrm{Area}(\Sigma)\leq\frac{16\pi(1+g)}{k}\displaystyle e^{\sup_{\Sigma}(f)-\inf_{\Sigma}(f)}=\frac{16\pi(1+g)}{k}e^{\textrm{osc}_{\Sigma}(f)}.
$$
That is, (\ref{area}) holds.

\qed

Now, suppose that $M$ is not simply connected.  We use a covering argument as in \cite{CS}.
\begin{prop} \label{est-area-1}
Let $(M^{3},\bar{g},e^{-f}d\mu)$ be a complete  smooth metric measure space  with $\overline{\Ric}_{f}\geq k>0$, where $k$ is a positive constant.  Let
$\Sigma^{2}$  be a closed embedded f-minimal surface. If $\Sigma$ is   contained in a  bounded domain $D$ of $M$ with  convex boundary $\partial D$,  then
\begin{equation}\label{f-area-1}
Area_f(\Sigma)\leq\frac{16\pi}{k}\biggr(\frac{2}{|\pi_{1}|}-\frac{1}{2}\chi(\Sigma)\biggr)e^{-\inf_\Sigma f}.
\end{equation}
\begin{equation}\label{area-1}
Area(\Sigma)\leq\frac{16\pi}{k}\biggr(\frac{2}{|\pi_{1}|}-\frac{1}{2}\chi(\Sigma)\biggr)e^{\text{osc}_\Sigma f}.
\end{equation}
where $|\pi_{1}|$ is the order of the first fundamental group of $M$, and $\chi(\Sigma)$ is the Euler characteristic of $\Sigma$.
\end{prop}
\noindent\textbf{Proof}. Let $\hat{M}$ be the universal covering manifoldof $M$. By Proposition \ref{group},  the covering is a finite $|\pi_{1}|$-fold covering.  Let $\hat{\Sigma}$ be the lifting of $\Sigma$.  In the proof of Theorem \ref{est1-1}, we have shown  that  $\hat{\Sigma}$ is orientable and   satisfies the assumption of Theorem \ref{est1-1}.  Hence
  Theorem~\ref{est1-1} implies that the first eigenvalue of the weighted Laplacian of $\hat{\Sigma}$ satisfies $\lambda_{1}({\hat\Delta}_{\hat{f}})\geq\frac{k}{2}$, where $\hat{f}$ is the lift of $f$. By Proposition \ref{est2}, we conclude that
\begin{equation*}
Area(\hat{\Sigma})\leq\frac{16\pi}{k}\biggr(2-\frac{1}{2}\chi(\hat{\Sigma})\biggr)e^{osc_{\hat{\Sigma}}(\tilde{f})},
\end{equation*}
and
\begin{equation*}
\textrm{Area}_{\hat{f}}(\hat{\Sigma})=\int_{\hat{\Sigma}}e^{-\hat{f}}d\sigma\leq\frac{16\pi}{k}\biggr(2-\frac{1}{2}\chi(\hat{\Sigma})\biggr)e^{-\inf_{\hat{\Sigma}}(\hat{f})}.
\end{equation*}
Thus (\ref{f-area-1}) and (\ref{area-1})  follow from the facts that $\chi(\hat{\Sigma})=|\pi_{1}|\cdot\chi(\Sigma)$, $\textrm{Area}(\hat{\Sigma})=|\pi_{1}|\cdot \textrm{Area}(\Sigma)$,  $\textrm{Area}_{\hat{f}}(\hat{\Sigma})=|\pi_{1}|\cdot \textrm{Area}_f(\Sigma)$,  $\textrm{inf}_{\hat{\Sigma}}(\hat{f})=\textrm{inf}_{\Sigma}(f)$ and $\textrm{osc}_{\hat{\Sigma}}(\hat{f})=\textrm{osc}_{\Sigma}(f)$.

\qed

In the following, we will give the upper bound for the total curvature of $f$-minimal surfaces. Here the term  {\it  the total curvature} of $\Sigma$ means $\int_{\Sigma}|A|^2d\sigma$ not $\int_{\Sigma}Kd\sigma$.

\begin{prop}\label{est-total}
If $(M^{3},\bar{g},e^{-f}d\mu)$ is a  smooth metric measure space with $\overline{\textrm{Ric}}_{f}\geq k$, where $k$ is a positive constant. Let   $\Sigma^{2}\subset M$ be  a closed embedded $f$-minimal surface with genus $g$.  If $\Sigma$ is   contained in a  bounded domain $D$ of $M$ with  convex boundary $\partial D$, then  $\Sigma$ satisfies
\begin{equation}
\int_{\Sigma}|A|^2d\sigma \leq C,
\end{equation}
where $A$ is the second fundamental form of $(\Sigma, \bar{g})$  and $C$ is a constant depending on the genus $g$ of $\Sigma$, the order $|\pi_1|$ of the first fundamental group of $M$, the maximum $\sup_{\Sigma}\bar{K}$ of the sectional curvature of $M$ on $\Sigma$, the lower bound $k$ of the Bakry-\'Emery Ricci curvature of $M$, the oscillation $\text{osc}_{\Sigma}(f)$ and the maximum $\sup_{\Sigma} |\bar{\nabla}f|$ on $\Sigma$.
\end{prop}
\noindent\textbf{Proof}. By the Gauss equation and Gauss-Bonnet formula,
\begin{align*}
\int_{\Sigma}|A|^2d\sigma&= \int_{\Sigma}H^2-2\int_{\Sigma}(K-\overline{K})\\
&= \int_{\Sigma}\langle \bar{\nabla}f, {\bf n}\rangle^2-4\pi\chi(\Sigma)+2\int_{\Sigma}\overline{K}\\
&\leq (\sup_{\Sigma} |\bar{\nabla}f|)^2\text{Area}(\Sigma)+8\pi(g-1)+2(\sup_{\Sigma}\bar{K})\text{Area}(\Sigma).
\end{align*}
Substituting   (\ref{area-1}) for $\text{Area}(\Sigma)$,  we have the conclusion of theorem.

\qed

To prove the compactness theorem in Section \ref{compactness1}, we need the following total curvature estimate for $(\Sigma, \tilde{g})$, which is a minimal  surface in $(M,\tilde{g})$.
\begin{prop}\label{est-total-conf}
If $(M^{3},\bar{g},e^{-f}d\mu)$ is a  smooth metric measure space with $\overline{\textrm{Ric}}_{f}\geq k$, where $k$ is a positive constant. Let  $\Sigma^{2}\subset M$ be a closed embedded $f$-minimal surface with genus $g$. If $\Sigma$ is contained in a  bounded domain $D$ of $M$ with convex boundary $\partial \Omega$,  then  $\Sigma$  satisfies that
\begin{equation}\label{total-conf}
\int_{\Sigma}|\tilde{A}|_{\tilde{g}}^2d\tilde{\sigma} \leq C,
\end{equation}
where $\tilde{A}$ is the second fundamental form of $(\Sigma,\tilde{g})$ with respect to the conformal metric $\tilde{g}=e^{-f}\bar{g}$ of $M$ and   C is a constant depending on the genus $g$ of $\Sigma$, the order $|\pi_1|$ of the first fundamental group of $M$, the maximum $\sup_{\Sigma}\tilde{K}$ of the sectional curvature of $(M,\tilde{g})$ on $\Sigma$, the lower bound $k$ of the Bakry-\'Emery Ricci curvature of $M$, the oscillation $\text{osc}_{\Sigma}(f)$ on $\Sigma$.
\end{prop}

\noindent\textbf{Proof}.  By the Gauss equation and Gauss-Bonnet formula, we have
\begin{align*}
\int_{\Sigma}|\tilde{A}|_{\tilde{g}}^2d\tilde{\sigma}&= \int_{\Sigma}\tilde{H}^2-2\int_{\Sigma}(\tilde{K}^{\Sigma}-\tilde{K}^M)d\tilde{\sigma}\\
&=-4\pi\chi(\Sigma)+2\int_{\Sigma}\tilde{K}d\tilde{\sigma}\\
&\leq 8\pi(g-1)+2(\sup_{\Sigma}\tilde{K})\text{Area}((\Sigma,\tilde{g}))\\
&=8\pi(g-1)+2(\sup_{\Sigma}\tilde{K})\text{Area}_f(\Sigma).
\end{align*}
In the above, we used $\tilde{H}=e^fH_f=0$ and $\text{Area}((\Sigma,\tilde{g}))=\text{Area}_f(\Sigma)$.
From (\ref{f-area-1}) in Proposition \ref{est-area-1},  Inequality (\ref{total-conf}) holds.

\qed

\section{Compactness of compact $f$-minimal surfaces}\label{compactness1}

We will prove  some compactness theorems for closed embedded $f$-minimal surfaces in $3$-manifolds. We have two ways to prove Theorem \ref{compact-bound}. 

The first  proof roughly follows  the one  in  \cite{CM4}  (cf  \cite{CS}) with some modifications. The modifications can be made because we have  the assumptions that $f$-minimal surfaces are contained in the closure  of a bounded domain ${\Omega}$  of $M$   and $\overline{\Omega}$ is contained in a  bounded domain $U$ with convex boundary. The second proof will need a compactness theorem of complete embedded $f$-minimal surfaces that was proved by us in \cite{CMZ}. 

We prefer to give two proofs here since the first one is independent of the compactness theorem of complete embedded $f$-minimal surfaces. But the compactness theorem of complete embedded $f$-minimal surfaces needs a theorem about non-existence of $L_f$-stable minimal surfaces (\cite{CMZ} Theorem 3).

{\it First proof}.

We first prove a  singular compactness theorem,  which is a variation of Choi-Schoen's  \cite{CS}   singular compactness theorem (cf  Proposition 7.14 in  \cite{CM4}, Anderson \cite{A} and White \cite{W}).  Namely,
\begin{prop} \label{singular} Let $(M^3,\bar{g})$ be a   $3$-manifold. Assume that $\Omega$  is  bounded  domain in $M$. Let $\Sigma_i $ be a sequence of closed embedded minimal surfaces  contained in  $\overline{\Omega}$, with  genus $g$,  and satisfying
\begin{equation}\label{area-0}\text{Area}(\Sigma_i)\leq C_1
\end{equation}
and 
\begin{equation}\label{total-1}\int_{\Sigma_i}|A_{\Sigma_i}|^2\leq C_2,
\end{equation}
Then there exists a finite set of points $\mathcal{S}\subset \bar{\Omega}$ and a subsequence, still denoted by  $\Sigma_{i}$, that converges uniformly in the $C^m$ topology ( $m\geq 2$ ) on compact subsets of $M\backslash\mathcal{S}$ to a complete minimal surface $\Sigma\subset \bar{\Omega}$ (possibly with multiplicity).  

The subsequence also converges to $\Sigma$ in extrinsic Hausdorff distance. $\Sigma$ is smooth, embedded in $M$, has genus at most $g$ and satisfies (\ref{area-0}) and (\ref{total-1}).
\end{prop}

\noindent \textbf{Proof.}  We may use the same argument as  that of Proposition 7.14 in  \cite{CM4}. Moreover,   $\Sigma_i\subset  \bar{\Omega}$ implies  that the singular set $S\subset \bar{\Omega}$  and the smooth surface $\Sigma\subset \bar{\Omega}$. Here we omit the details of  proof.

\qed

We  can apply Proposition \ref{singular} to the $f$-minimal surfaces which are minimal in the conformal metric.
\begin{lemma} \label{sing-conf} Let $(M^{3},\bar{g},e^{-f}d\mu)$ be a  smooth  metric measure space.   Assume that $\Omega$ is a  bounded  domain in $M$. Let $\Sigma_i\subset \overline{\Omega} $ be a sequence of closed embedded $f$-minimal surfaces of genus $g$. Suppose that $\tilde{g}=e^{-f}\bar{g}$ on $M$ and $(\Sigma_i, \tilde{g})$ satisfy that
\begin{equation}
\text{Area}\label{area-2}((\Sigma_i,\tilde{g}))=\text{Area}_f(\Sigma_i)\leq C_1
\end{equation}
and 
\begin{equation}\label{total-2}\int_{\Sigma_i}|\tilde{A}_{\Sigma_i}|_{\tilde{g}}^2d\tilde{\sigma}\leq C_2,
\end{equation}
where $\tilde{A}_{\Sigma_i}$ and $d\tilde{\sigma}$ denote the second fundamental form and the volume element of $(\Sigma_i,\tilde{g})$ respectively.
Then there exists a finite set of points $\mathcal{S}\subset \bar{\Omega}$ and a subsequence, still denoted by  $\Sigma_{i}$, that converges uniformly in the $C^m$ topology ( $m\geq 2$) on compact subsets of $M\backslash\mathcal{S}$ to a complete  $f$-minimal surface $\Sigma\subset \overline{\Omega}$   (possibly with multiplicity). 

 The subsequence also converges to $\Sigma$ in extrinsic Hausdorff distance. $\Sigma$ is smooth, embedded in $M$, has  genus at most $g$, and satisfies (\ref{area-2}) and (\ref{total-2}).
\end{lemma}

\noindent\textbf{Proof}.  Since an $f$-minimal surface in the original metric $\bar{g}$ is equivalent to it is minimal in the conformal metric $\tilde{g}$, we can apply  Proposition \ref{singular} to get the conclusion of the lemma.

\qed

\noindent\textbf{Proof of Theorem \ref{compact-bound}}.  We first consider that $M$ is simply connected. Since $\Sigma_i\subset \overline{\Omega}\subset  U$. By Proposition  \ref{est2} and Proposition  \ref{est-total-conf},
$$\text{Area}((\Sigma_i,\tilde{g}))=\text{Area}_f(\Sigma_i)\leq C_1$$
and $$\int_{\Sigma_i}|\tilde{A}_{\Sigma_i}|_{\tilde{g}}^2d\sigma_{\tilde{g}}\leq C_2,$$
where $C_1$ and $C_2$  depend on $g$, $\sup_{\Omega_j}f$, $\sup_{\Omega_j}\tilde{K}$ and $k$.

By Lemma \ref{sing-conf},   there exists a finite set of points $\mathcal{S}\subset \tilde{\Omega}$ and a subsequence $\Sigma_{i'}$ that converges uniformly in the $C^m$ topology (any $m\geq 2$) on compact subsets of $M\backslash\mathcal{S}$ to a complete  $f$-minimal surface $\Sigma\subset \bar{\Omega}$ without boundary (possibly with multiplicity).  $\Sigma$ is smooth, embedded  in $M$ and  has genus at most $g$.  Equivalently, with respect to the conformal metric $\tilde{g}$, a subsequence $\Sigma_{i'}$  of minimal surffaces  converges uniformly in the $C^m$ topology  on compact subsets of $M\backslash\mathcal{S}$ to a complete minimal surface $\Sigma$, where $\Sigma\subset \bar{\Omega}$.

Since complete embedded $\Sigma\subset\overline{\Omega}$ satisfies (\ref {area-2}),  it must be properly embedded (Proposition \ref{prop-imer1}), thus 
closed, and orientable.

We need to  prove that the convergence is
smooth across the points $\mathcal{S}$. By Allard's regularity theorem, it suffices to prove that the convergence has multiplicity one. If the multiplicity is not one, by a proof similar to that of Choi-Schoen \cite{CS}  (also cf \cite{CM4} P. 249), we can show that there is an $i$ big enough and a $\Sigma_i$ in the convergent subsequence, so that the first eigenvalue of the Laplacian $\tilde{\Delta}^{\Sigma_i}$ on $\Sigma_i$ with the conformal metric $\tilde{g}$ satisfies $\lambda_1(\tilde{\Delta}^{\Sigma_i})<\frac k2e^{\inf_\Omega f}$.   We have
\begin{align*}\lambda_1(\tilde{\Delta}^{\Sigma_i})&=\displaystyle\inf\{\frac{\int_{\Sigma_i}|\tilde{\nabla} \phi |^2_{\tilde{g}}d\tilde{\sigma}}{\int_{\Sigma_i}\phi^2d\tilde{\sigma}}, \int_{\Sigma_i}\phi d\tilde{\sigma}=0\}\\
&=\displaystyle\inf\{\frac{\int_{\Sigma_i}|\nabla\phi |^2d\sigma}{\int_{\Sigma_i}\phi^2e^{-f}d\sigma}, \int_{\Sigma_i}\phi e^{-f}d\sigma=0\}\\
&\geq \lambda_1(\Delta_f^{\Sigma_i})e^{\inf_\Omega f}.
\end{align*}

By Theorem \ref{est1-1},  $\Sigma_i\subset\bar \Omega\subset U$ implies $ \lambda_1(\Delta_f^{\Sigma_i} )\geq \frac k2.$ Thus we have
a contradiction.

When  $M$ is not simply connected, we use a covering argument. The assumption of $\overline{\text{Ric}}_f\geq k$, where constant $k>0$, implies that $M$ has finite  fundamental group $\pi_1$ (Proposition \ref{group}). We consider the finite-fold universal covering $\hat{M}$. By the proof of Theorem \ref{est1-1}, we know that the corresponding lifts of $\Sigma_i$, $\overline{\Omega}$ and $U$ satisfy that $\hat{\Sigma}_i\subset\hat{\overline{\Omega}}\subset\hat{U}$.  Then Propositions \ref{est-area-1} and  \ref{est-total-conf} give the uniform bounds of area and total curvature in the conformal metric $\hat{\tilde{g}}$ on $\hat{M}$. 
By the assertion on the simply connected ambient manifold before,  we have the smooth convergence of a subsequence of $\hat{\Sigma}_i$. This implies the smooth convergence of a subsequence of $\Sigma_i$.

\qed

{\it Second Proof}. In \cite{CMZ}, we proved the following
\begin{theorem}\label{prin-intro}\cite{CMZ}
Let $(M^{3},\bar{g},e^{-f}d\mu)$ be a complete smooth metric measure space and $\overline{\textrm{Ric}}_{f}\geq k$, where $k$ is a positive constant. Given an integer $g\geq0$ and a constant $V>0$, the space $S_{g, V}$ of smooth complete embedded $f$-minimal surfaces $\Sigma\subset M$ with 
\begin{itemize}
\item genus at most $g$,
\item $\partial\Sigma=\emptyset$, 
\item  $\int_{\Sigma}e^{-f}d\sigma\leq V$
\end{itemize}
is compact in the $C^m$ topology, for any $m\geq 2$.
Namely, any sequence of $S_{g, V}$ has a subsequence that converges in the $C^m$ topology on compact subsets to a surface in $S_{D, g}$, for any $m\geq 2$.
\end{theorem}

\noindent \textbf{Proof of Theorem \ref{compact-bound}}.  Since a surface in $S$ is contained in $\bar{\Omega}\subset U$, by Proposition \ref{est-area-1}, we have the uniform bound $V$ of the weighted volume of closed embedded $f$-minimal surfaces in $S$. Hence Theorem \ref{prin-intro} can be applied.  Moreover $\Sigma_i\subset\overline{\Omega}$ implies that the smooth limit surface $\Sigma\subset\overline{\Omega}$. Otherwise, 
 since  the subsequence $\{\Sigma_i\}$ converges uniformly in the $C^m$ topology ($m\geq 2$)  on any compact subset of $M$ to $\Sigma$, there is a  surface $\Sigma_i$  (with  index $i$ big enough) in the subsequence would not satisfy $\Sigma_i\subset\overline{\Omega}$.  
 
 By Proposition \ref{prop-imer1}, $\Sigma$ must be properly embedded.  Thus $\Sigma$ must be closed.
 
 \qed

Using Theorem \ref{compact-bound}, we may prove Theorem \ref{compact-intro} in Introduction.

\begin{lemma} \label{bound}Let $(M^{3},\bar{g},e^{-f}d\mu)$ be a complete noncompact  smooth  metric measure space with $\overline{\textrm{Ric}}_{f}\geq k>0$. If $\Sigma$ is any closed  $f$-minimal surface in $M$  with genus  at most $g$ and diameter at most $D$,  then $\Sigma\subset \overline{B}_r(p )$ for some $r>0$ (indepedent of $\Sigma$), where $B_r( p)$ is a ball in $M$ with radius $r$ centered at $p\in M$.
\end{lemma}
\noindent\textbf{Proof}. Fix a closed $f$-minimal surface $\Sigma_0$. Obviously, $\Sigma_0\subset B_{r_0}(p )$ for some $r_0>0$. Proposition \ref{inter} says that  $\Sigma$ and $\Sigma_0$ must intersect. Then for $x\in \Sigma$, $$d(p, x)\leq d(p, x_0)+d(x_0, x)\leq r_0+D, x_0\in \Sigma_0.$$
Taking  $r=r_0+D$, we have $\Sigma\subset \overline{B}_{r_0+D}.$

\qed

\begin{remark}\label{extrinsic}  In Lemma \ref{bound} and hence in Theorem \ref{compact-intro}, $D$ is a bound of intrinsic  diameter of closed $f$-minimal surfaces or  a bound of extrinsic  diameter of closed $f$-minimal surfaces.  Also by Proposition \ref{inter}, the assumption that $f$-minimal surfaces are contained in the closure of a bounded domain $\Omega$ in Theorem \ref{compact-bound}  is equivalent to  that the uniform upper bound of the extrinsic  diameter of $f$-minimal surfaces.
\end{remark}

\noindent\textbf{Proof of Theorem \ref{compact-intro}}. By Lemma \ref{bound}, we may apply Theorem \ref{compact-bound} to the space $S_{D,g}$. Next, 
the closed embedded limit $\Sigma$ must have diameter at most $D$. Otherwise, since  the subsequence $\{\Sigma_i\}$ converges uniformly in the $C^m$ topology ($m\geq 2$)  on any compact subset of $M$ to $\Sigma$, there is a  surface $\Sigma_i$  (with the index $i$ big enough) in the subsequence would have diameter greater than $D$. So $\Sigma$ must in $S_{D,g}$

\qed

\section{Appendix}

In this section,  we include the proof Proposition \ref{prop1} in Section \ref{reilly}.

\noindent\textbf{Proof of Proposition \ref{prop1}}.  Recall the Bochner formula
\begin{eqnarray*}
\frac{1}{2}\bar{\Delta}_{f}|\overline{\nabla}u|^{2}-\langle\overline{\nabla}u,\overline{\nabla}(\bar{\Delta}_{f}u)\rangle=|\overline{\nabla}^{2}u|^{2}+
\overline{\textrm{Ric}}_{f}(\overline{\nabla}u,\overline{\nabla}u).
\end{eqnarray*}
Integrating this equation on $\Omega$ with respect to  weighted measure $e^{-f}d\mu$, we obtain
\begin{equation*}
\int_{\Omega}\bigr(\frac{1}{2}\bar{\triangle}_{f}|\overline{\nabla}u|^{2}-\langle\overline{\nabla}u,\overline{\nabla}(\bar{\triangle}_{f}u)\rangle\bigr)e^{-f}=
\int_{\Omega}|\overline{\nabla}^{2}u|^{2}e^{-f}+\int_{\Omega}\overline{\textrm{Ric}}_{f}(\overline{\nabla}u,\overline{\nabla}u)e^{-f}.\\
\end{equation*}
On the other hand,  by the divergence formula, it holds that
\begin{align*}
&  \frac{1}{2}\bar{\Delta}_{f}|\overline{\nabla}u|^{2}-\langle\overline{\nabla}u,\overline{\nabla}(\bar{\Delta}_{f}u)\rangle\\
&\quad=\frac{1}{2}\overline{\textrm{div}}\bigr(e^{-f}
\overline{\nabla}|\overline{\nabla}u|^{2}\bigr)e^{f}-\overline{\textrm{div}}\bigr(e^{-f}\bar{\Delta}_{f}(u)\overline{\nabla}u\bigr)e^{f}+(\bar{\Delta}_{f}u)^{2}. 
\end{align*}
Integrating and applying the Stokes'  theorem,  we have
\begin{align}\label{append-1}
&\int_{\Omega}\bigr(\frac{1}{2}\bar{\Delta}_{f}|\overline{\nabla}u|^{2}-\langle\overline{\nabla}u,\overline{\nabla}(\bar{\Delta}_{f}u)\rangle\bigr)e^{-f}\\&\quad=
\int_{\partial\Omega}\bigr(\frac{1}{2}|\overline{\nabla}u|^{2}_{\nu}-(\bar{\Delta}_{f}u)u_{\nu}\bigr)e^{-f}
+\int_{\Omega}(\bar{\Delta}_{f}u)^{2}e^{-f}.\nonumber
\end{align}
Then
\begin{align}\label{append-2}
\frac{1}{2}|\overline{\nabla}u|^{2}_{\nu}-(\bar{\Delta}_{f}u)u_{\nu}
&=\langle\overline{\nabla}_{\nu}\overline{\nabla}u,\overline{\nabla}u\rangle-(\bar{\Delta}_{f}u)u_{\nu}\\
&=\langle\overline{\nabla}_{\overline{\nabla}u}\overline{\nabla}u,\nu\rangle-(\bar{\Delta}_{f}u)u_{\nu}\nonumber\\
&=\langle\overline{\nabla}_{\nu}\overline{\nabla}u,\nu\rangle u_{\nu}+\langle\overline{\nabla}_{\nabla u}\overline{\nabla}u,\nu\rangle-(\bar{\Delta}_{f}u)u_{\nu}\nonumber\\
&=\bigr(\langle\overline{\nabla}_{\nu}\overline{\nabla}u,\nu\rangle-\bar{\Delta}u+\langle\overline{\nabla}f,\overline{\nabla}u\rangle\bigr) u_{\nu}\nonumber\\
&\quad \quad+\langle\nabla u,\nabla u_{\nu}\rangle-
\langle\overline{\nabla}u,\overline{\nabla}_{\nabla u}\nu\rangle\nonumber\\
&=\bigr(-\Delta u-Hu_{\nu}+\langle\nabla f,\nabla u\rangle+\langle\bar{\nabla}f,\nu\rangle u_{\nu}\bigr)u_{\nu}\nonumber\\
&\quad\quad +\langle\nabla u,\nabla u_{\nu}\rangle-\langle\nabla u,\overline{\nabla}_{\nabla u}\nu\rangle\nonumber\\
&=-(\Delta_{f}u+H_{f}u_{\nu})u_{\nu}+\langle\nabla u,\nabla u_{\nu}\rangle-A(\nabla u,\nabla u).\nonumber
\end{align}
where $H_{f}=H-\langle\overline{\nabla}f,\nu\rangle$.
By substituting (\ref{append-2}) into (\ref{append-2}), we obtain
\begin{align*}
&\int_{\Omega}\bigr(\frac{1}{2}\bar{\Delta}_{f}|\overline{\nabla}u|^{2}-\langle\overline{\nabla}u,\overline{\nabla}(\bar{\Delta}_{f}u)\rangle\bigr)e^{-f}\\
&=-\int_{\partial\Omega}(\Delta_{f}u)u_{\nu}e^{-f}-\int_{\partial\Omega}H_{f}u_{\nu}^{2}e^{-f} \\
&\quad\quad+\int_{\partial\Omega}(\langle\nabla u,\nabla u_{\nu}\rangle-A(\nabla u,\nabla u))e^{-f} +\int_{\Omega}(\bar{\Delta}_{f}u)^{2}e^{-f}\\
&=-2\int_{\partial\Omega}(\Delta_{f}u)u_{\nu}e^{-f}-\int_{\partial\Omega}H_{f}u_{\nu}^{2}e^{-f}\\
&\quad\quad -\int_{\partial\Omega}\hspace{-2mm}A(\nabla u,\nabla u)e^{-f}+\int_{\Omega}(\bar{\Delta}_{f}u)^{2}e^{-f}.
\end{align*}
Consequently
\begin{eqnarray*}
\int_{\Omega}(\bar{\Delta}_{f}u)^{2}e^{-f}&=&\int_{\Omega}|\overline{\nabla}^{2}u|^{2}e^{-f}+\int_{\Omega}\overline{\textrm{Ric}}_{f}(\overline{\nabla}u
,\overline{\nabla}u)e^{-f}+2\int_{\partial\Omega}u_{\nu}\Delta_{f}(u)e^{-f}\\
& &+\int_{\partial\Omega}A(\nabla u,\nabla u)e^{-f}+\int_{\partial\Omega}u_{\nu}^{2}H_{f}e^{-f}.
\end{eqnarray*}
\qed

\end{document}